\documentclass[aps,apl,reprint,groupedaddress, showkeys,showpacs]{revtex4-1}
\usepackage[utf8]{inputenc}
\usepackage{amsmath,amssymb}
\usepackage{graphicx,float}

\begin{document}

\title{Emergence of the giant weak-component in directed random graphs with arbitrary degree distributions}
\author{Ivan Kryven} 
\email{i.kryven@uva.nl}
\affiliation{University of Amsterdam, PO box  94214, 1090 GE, Amsterdam, The Netherlands} 

\begin{abstract}
The weak component generalizes the idea of connected components to directed graphs. In this paper, an exact criterion for existence of the giant weak component is derived for directed graphs with arbitrary bivariate degree distributions.
In addition we consider a random process for evolving directed graphs with bounded degrees. 
The bounds are not the same for different vertices but satisfy a pre-defined distribution.
The analytic expression obtained for the evolving degree distribution is then combined with the weak-component criterion to obtain the exact time of the phase transition. The phase-transition time is obtained as a function of the distribution that bounds the degrees.
Remarkably, when viewed from the step polymerization formalism, the new results yield Flory-Stockmayer gelation theory and generalize it to a broader scope.

\end{abstract}


\maketitle

\section{Introduction}
From reactions fueling cells in our bodies to internet links binding the World Wide Web into a small-world structure, networks are at the basis of many phenomena.
Random Graph theory sets up a common toolbox studying networks independently of their context from a probabilistic point of view. The most basic random graph model was introduced by Erd\H{o}s\cite{Erdos}. This model refers to a set of vertices, where the probability of two vertices being connected is chosen in advance as the only model parameter. Vertices from such a graph satisfy a specific degree distribution, namely the Poisson distribution. Since Erd\H{o}s' first results appeared, many models yielding other degree distributions followed (see for example \cite{bollobas1998} and the citations therein). This exploration for new models was driven by both a theoretical desire and a practical necessity.

With respect to soft-matter physics, among other disciplines, the random graph theory has been a source of inspiration for coagulation, polymerization, and elasticity models. 
In this frameworks the set of all admissible graphs yielded by a model is drastically constrained by the physical context.
There are three main tools providing a means to impose these constrains:  bond percolation on (pseudo) lattices (\textit{e.g.} as in the study on Bethe lattice by Fisher and Essam \cite{fisher1961}), kinetic perspective on random graphs, and imposing a pre-defined degree distribution.

Some notable examples of the kinetic perspective on random graphs include, but are not limited to: analytic results obtained by Ben-Naim \cite{ben2005,ben2011,ben2011a}, Lushnikov \cite{leyvraz2015,lushnikov2015,lushnikov2015a}, Buffet\cite{Buffet1991}, Gordon \cite{gordon1962} and their coauthors; numerical studies of step-growth and cross-linking polymerization by Kryven et al.\cite{kryven2014a,kryven2014c,kryven2013c},  and algorithmic method introduced by Hillegers and Slot \cite{hillegers2016,hillegers2015}. In the above enlisted cases  (with exception of Refs. \cite{ben2011,hillegers2015}) the degree distribution is either very simple (0, 1 or 2 edges) or not constrained at all.
Models that implement restrictions on the degree distribution in the evolving undirected random graphs are considered in \cite{rucinski1992,ben2011}. This restriction is  implemented as a single upper bound on the degrees.
In the same time, the algorithmic study on directed graphs (inspired by the polymerization structures) \cite{hillegers2015} does allow one to impose the degree bounds as a distribution, yet the algorithm is applicable only prior to the phase transition. 

Clearly, a degree distribution does not define a graph uniquely. That said, an attractive alternative to the classical models is to define a random graph by a given degree distribution assuming that apart from the degree distribution the graph is absolutely random. 
This line of research was introduced by Molloy and Reed\cite{molloy1995} and was later developed further by Newman, Strogatz and Watts\cite{newman2001}.
Studying properties of random graphs defined by their degree distribution is not simply an abstract problem, it has a clear practical motivation. For instance, one may consider an empirical degree distribution that is based on measured or observed data. An observer collecting such data is likely to be either embedded into the network himself, thus viewing it locally, or to be distanced far apart, thus observing only the global properties. Indeed, one may study individual servers of the Internet but the question of the global connectivity structure is far less trivial \cite{broder2000}, or one may observe global properties of a complex polymer material without exhaustive knowledge on how the individual molecules are interconnected. Expressing global properties of random graphs in terms of their degree distribution builds up an essential link between the local and the global.

In undirected graphs a connected component is a set of all vertices that can be reached from a given vertex by following the edges recursively. Many random graphs are known to experience the phase transition: the point when a connected component that has size of the same order of magnitude as the whole network emerges (the giant component).
The idea of the connected component can also be generalized to directed graphs, \textit{i.e.} graphs having all edges with a specific direction. For a selected vertex:
(a)\textit{out-component} is a set of vertices that can be reached by recursively following all out-edges forward; 
(b)\textit{in-component} is a set of vertices that can be reached by recursively following all in-edges backwards; 
(c)\textit{weak component} is a set of vertices that can be reached by recursively following all edges regardless their orientation.

Even when focusing on weak components alone, one finds many applied studies exploiting the concept; for instance in epidemiology\cite{firestone2012,eubank2004}, data mining\cite{gunnemann2010,kirsch2006}, communication networks\cite{butts2007}  exploring World Wide Web structure\cite{broder2000,dorogovtsev2001}, \textit{etc}.
So what information on component sizes one may obtain just by knowing the degree distribution of a directed random graph?
In the case of undirected graphs the question has been answered by Molloy and Reed\cite{molloy1998}.  
A link between the degree distribution and some properties of directed graphs has been showen by Chen et al.\cite{chen2013}. 
A new theory studying sizes of in- and out-components was introduced by Newman et al.\cite{newman2001}. 
A connection between the degree distribution and the giant weak component, however, has not been investigated in depth.
Moreover, some authors (\textit{e.g.} in \cite{newman2001,dorogovtsev2001}) intentionally do not study weak components separately, arguing that in this case, the graph effectively becomes undirected and should be treated with the known formalism.  
This statement, generally speaking, is a misconception since even though we disregard directional information when calculating size of the weak component, direction of the edges does affect the topology of the network. 

The current paper is organized in two parts. First, a correct criterion for existence of the giant weak component will be derived. This criterion shapes as an inequality involving moments of the degree distribution. This result complements the prior findings on in-/out-components and components in undirected graphs. The criterion can be immediately applied to simulated or empirical degree distributions.
In the second part of the paper, an analytic expression for the bivariate degree distribution is derived for a specific random-graph time process. 
In this process the directed random graph evolves starting from a set of disconnected vertices.
Similarly to undirected case considered in \cite{ben2011}, the degree of a vertex is bounded, but the bounds are not the same for different vertices. Therefore, we deal with a bivariate distribution of bounds as an input parameter.
The probability of a vertex to receive an edge is proportional to the difference between the bound and the actual number of   edges that are incident to the vertex.
 The weak-component criterion is then applied to obtain the phase-transition time as a function of the input parameters. Remarkably, these results produce the Flory-Stockmayer gelation theory \cite{flory1941,Stockmayer1944,ziff1980} as a special case and thus constitute a more general theory for gelation in themselves.

\section{Criterion of the phase transition for an arbitrary degree distribution}
In an undirected graph, the degree distribution defines the probability of having a specific number of edges for a randomly selected vertex.
In a directed graph, each vertex has an in-degree and out-degree that count edges coming to and leaving from the vertex.
For a given directed random graph, a bivariate degree distribution, $u(n,k),\; n,k=0,1,\dots$ denotes the probability that a randomly chosen vertex has in-degree $n$ and out-degree $k$. There are two extra properties that $u(n,k)$ has to satisfy to be a valid degree distribution. The total probability has to sum up to unity, 
$\sum_{n,k=1}^{\infty} u(n,k)=1,$ and the total numbers of in-edges and out-edges have to coincide, $\sum_{n,k=1}^{\infty}(n-k) u(n,k)=0.$
Let $\mu_{ij}$ denote partial moments of $u(n,k),$ 
\begin{equation}
\label{eq:MoM}
\mu_{ij} = \sum_{n,k}n^i k^j u(n,k).
\end{equation}
The edge balance  and the normalization condition for $u(n,k)$ can be rewritten using the moment notation,
\begin{equation}
\label{eq:requirement}
\mu_{00}=1 \;\text{and}\; \mu_{10}=\mu_{01}=\mu.
\end{equation}
 It is often more convenient to work with generating functions than actual distributions. The  generating function for the degree distribution is defined as
\begin{equation}
\label{eq:transform}
U(z,w) =  \sum_{n=0,k=0}^\infty u(n,k) z^n w^k, \; z,w\in\mathbb{C},
\end{equation}
where $|z| \leq 1,|w|\leq 1$.
Alternatively, one may rewrite the equalities \eqref{eq:requirement} in terms of the generating function using combination of differentiation and evaluation at point (1,1), 
$$ U(z,w)\Big|_{z=1,w=1}=1, $$
and 
$$ \left(\frac{\partial}{\partial z}\textendash\frac{\partial}{\partial w} \right)U(z,w)\Big|_{z=1,w=1}=0.$$

\begin{figure}
\includegraphics[width=0.48\textwidth]{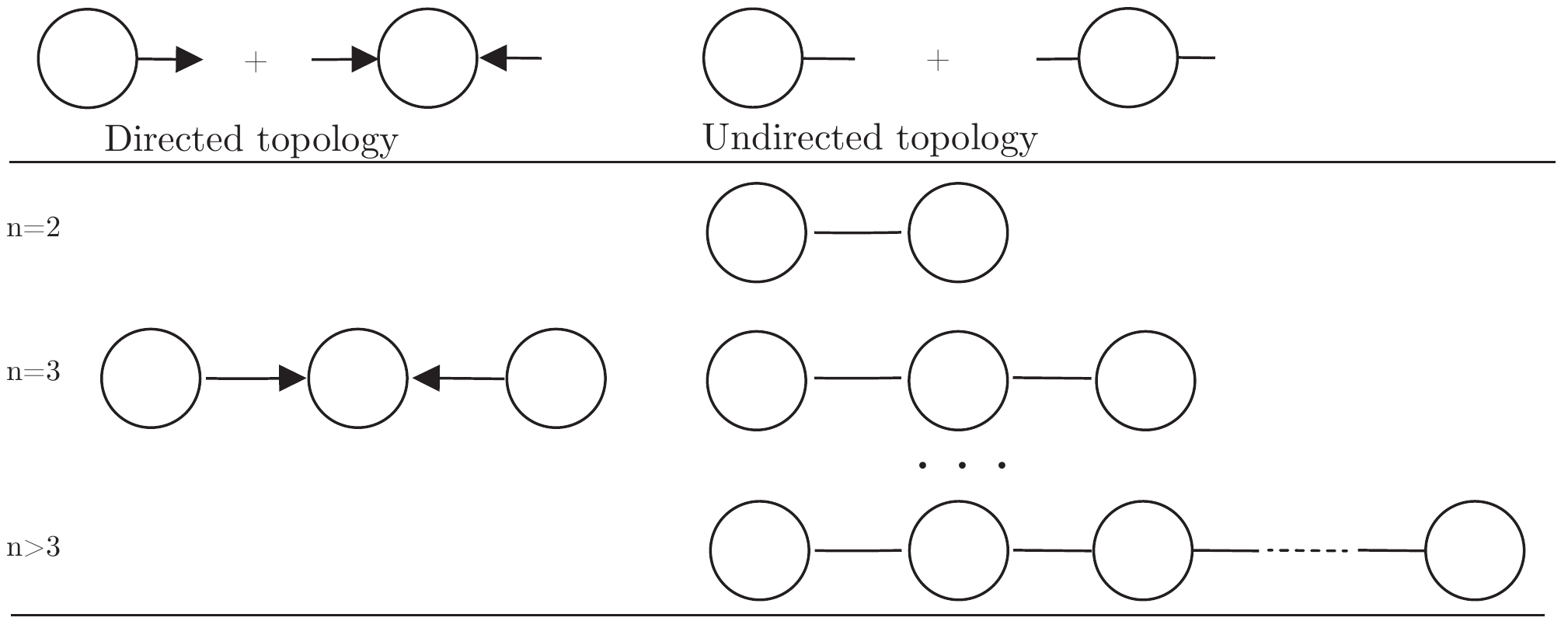}
\caption{An example illustrating how topologies of a directed graph and the undirected one satisfying an induced degree distribution, may be  drastically different. 
The directed graph on the left consists of vertices that have either one out-edge or two in-edges. 
The corresponding undirected graph on the right consist of vertices that have either one edge or two.
}
\label{fig:ex1}
\end{figure}

Now, let us introduce a bias into the process of vertex selection. Suppose, we select a vertex that is at the end of a randomly-chosen edge. The degree of the vertex is no longer governed by $u(n,k)$  since vertices of higher in-degree are more likely to be sampled. The correct degree distribution in this case is $u_{\text{in}}(n,k)=\frac{n}{\mu} u(n,k),$ which is generated by
\begin{equation}
\label{eq:Ux}
U_{\text{in}}(z,w)=\mu^{-1}\frac{\partial}{\partial z}U(z,w).
\end{equation}
In similar fashion, consider selecting a vertex that is at the beginning of a randomly-chosen edge. The degree distribution for such vertices is given by $u_{\text{out}}(n,k)=\frac{k}{\mu} u(n,k),$ which is generated by  
\begin{equation}
\label{eq:Uy}
U_{\text{out}}(z,w)=\mu^{-1}\frac{\partial}{\partial w}U(z,w).
\end{equation}
The weak component is a set of vertices that can be reached by recursively following all edges regardless their orientation.
When a  directed random graph is defined by the degree distribution only, the temptation is to say that the distribution of weak-component sizes is essentially the same as the distribution of component sizes that corresponds to undirected degree distribution, $d(l)=\sum_{n+k=l} u(n,k),\;l\in \mathbb{N}_0.$
This statement, generally speaking, is not correct since even though we disregard the directional information when calculating size of the weak component, direction of the edges does affect the topology of the network. This fact can be illustrated on a simple example: consider a bivariate degree distribution that is zero everywhere except for $u(1,0)=\frac{2}{3},\; u(0,2)=\frac{1}{3}.$ The directed random graph generated by such a distribution has only components of size 3. On another hand, component sizes in the undirected graph generated by the corresponding degree distribution ($d(1)=\frac{2}{3},\;d(2)=\frac{1}{3})$ are not bounded at all, see Fig.~\ref{fig:ex1}.

We will now extend the approach presented in \cite{newman2001} to cover the case of weak components for the directed graphs.
For a randomly selected vertex, let $w(n),\; n \in \mathbb{N}_0,$ such that $\sum_n w(n)=1,$ denotes the distribution of weak-component sizes; $w(n)$ is generated by $W(z)$. 
Analogously to definition of the distributions \eqref{eq:Ux},\eqref{eq:Uy}, consider a biased choice for the starting vertex. 
Suppose, one chooses an edge at random, and then selects the terminal vertex of this edge as a root.  In this case, let $w_{\text{in}}(n)$ (generated by $W_{\text{in}}(z)$) denotes the distribution of weak-component sizes associated with the root.
As another extreme, suppose one chooses an edge at random and then selects the source vertex of this edge as a root.
Similarly to the prior case,  let $w_{\text{out}}(n)$ (generated by $W_{\text{out}}(z)$) denotes the distribution of weak-component sizes associated with the root.
  \begin{figure}[htbp]
\begin{center}
\includegraphics[width=0.48\textwidth]{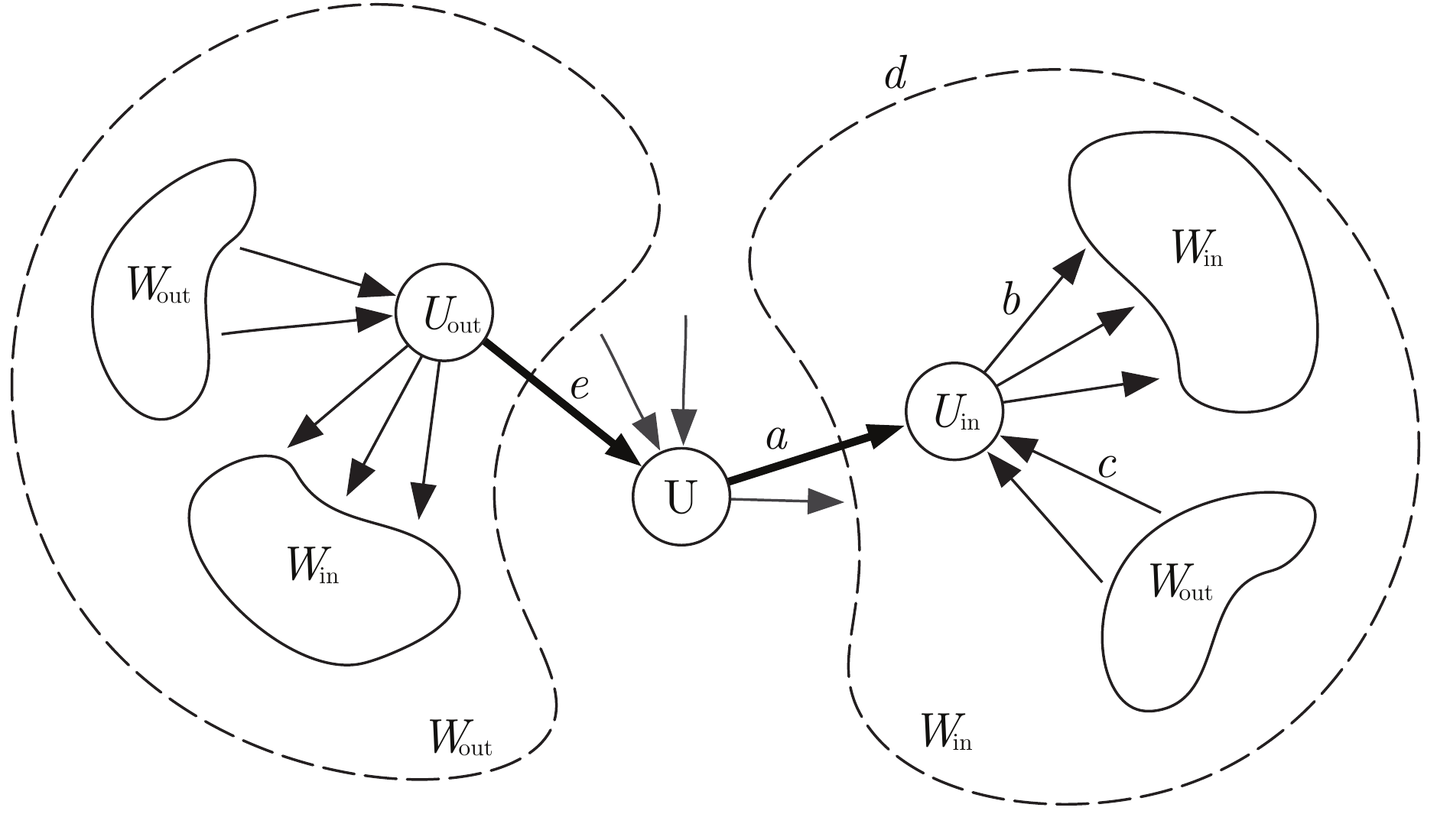}
\caption{Construction of the implicit equation for $w(n),$ the function generating weak-component sizes.}
\label{fig:clouds}
\end{center}
\end{figure}
The next step is to derive equalities binding $W_{\text{in}}, \, W_{\text{out}},\,U_{\text{in}},\, U_{\text{out}}$ together. 
Let us start by selecting a vertex (root) that we arrive at by following a random edge (edge $a$ in Fig.~\ref{fig:clouds}). 
According to the definition \eqref{eq:Ux} the probability of the root to have $n$ in-edges and $k$ out-edges is $u_{\text{in}}(n,k).$ 
 Each of the out-edges leaving the root is associated with a weak component of the size $w_{\text{in}}(n)$ (edges $b$ in Fig.~\ref{fig:clouds}), thus the sum of sizes of all components reached through the out-edges is distributed according to $k$-fold convolution $w_{\text{in}}(n) *w_{\text{in}}(n)*\dots* w_{\text{in}}(n).$ This sum is generated by $W_{\text{in}}(z)^k.$ A similar argument is constructed for all in-edges (edges $c$ in Fig.~\ref{fig:clouds}): the sum of sizes for all components reached through the in-edges is generated by $W_{\text{out}}(z)^n$. 
 A branch of such an exploration process will terminate in one of the two cases: a)when a vertex with at least one in-edge and no out-edges is reached (this happens with probability $U_{\text{in}}(1,0)$); b)when a vertex with at least one out-edge and no in-edges is reached (probability  $U_{\text{out}}(0,1)).$
 The distribution for the sum of sizes of all components originated at the root (\textit{i.e.} being reached through either in- or out-edge) is obtained as a summation over all possible configurations $(n,k),$ 
$$  \sum_{n,k} u_{\text{in}}(n,k) W_{\text{out}}(z)^n W_{\text{in}}(z)^k. $$
 Interestingly, this summation in itself can be viewed as a bivariate generating function of the type \eqref{eq:Ux} evaluated at point $z=W_{\text{out}}(z),\; w=W_{\text{in}}(z),$
 $$
 \sum_{n,k} u_{\text{in}}(n,k) W_{\text{out}}(z)^n W_{\text{in}}(z)^k = U_{\text{in}}\Big( W_{\text{out}}(z), W_{\text{in}}(z) \Big).
 $$
On another hand, the total number of all vertices reachable from the root plus one (component $d$ in Fig.~\ref{fig:clouds}) can also be considered as the size of the weak component reached by following an edge forwards. Thus one obtains a recurrence relation,
\begin{equation}
\label{eq:Win}
W_{\text{in}}(z) = z U_{\text{in}} \Big( W_{\text{out}}(z), W_{\text{in}}(z) \Big),
\end{equation}
where factor $z$  provides a unit translation in the component size distribution in order to  include the root into the component itself. 
A similar argumentation holds for $W_{\text{out}}(z)$.  Suppose one selects an edge at random and follows it in reverse (edge $e$ in Fig.~\ref{fig:clouds}) to reach a new root vertex. The degree of the root is described by $u_{\text{out}}(n,k)$ The sum of sizes for the weak components reached by the out-edges is generated by $W_{out}(z)^k$ and for in-edges this number is generated by $W_{in}(z)^n.$  The size of the whole weak component associated with the root is,
$$ \sum_{n,k} u_{\text{out}}(n,k) W_{\text{out}}(z)^n W_{\text{in}}(z)^k = U_{\text{out}}\Big( W_{\text{out}}(z), W_{\text{in}}(z) \Big).$$
Translating this distribution by unity yields the generating function for sizes of weak components that are reached by following an edge backwards,
\begin{equation}
\label{eq:Wout}
W_{\text{out}}(z) = z U_{\text{out}}\Big( W_{\text{out}}(z), W_{\text{in}}(z) \Big).
\end{equation}
When combined, Equations \eqref{eq:Win} and \eqref{eq:Wout} provide a sufficient means to uniquely define generating functions $W_{\text{out}}(z), W_{\text{in}}(z).$
Finally, we transit from sizes of biased weak components to sizes of weak components, generated by function $W(z).$ Consider a randomly selected vertex. Its degree distribution is generated by $U(z,w).$ The total sum of all component sizes reached via in- and out-edges plus 1 is generated by 
\begin{equation}
\label{eq:W}
 W(z) = z U\Big( W_{\text{out}}(z), W_{\text{in}}(z) \Big),
 \end{equation}
 which is the generating function for the weak-component size distribution. Similar relation for the giant in-component was derived in \cite{newman2001}.
 Even though, the triple \eqref{eq:Win},\eqref{eq:Wout},\eqref{eq:W} defines $W(z)$ implicitly, some properties of $W(z)$ may be extracted in an explicit form. For instance, we may find out if the random graph contains the giant weak component.
 
 Recalling that $U(1)=U_{\text{in}}(1)=U_{\text{out}}(1)=W(1)=1,$ the average size of the weak component to which a randomly chosen vertex belongs is given by,
\begin{equation}
 \begin{aligned}
 \label{eq:derW}
 W'(1)= &\left(z U\big( W_{\text{out}}(z), W_{\text{in}}(z) \big) \right)'\Big|_{z=1}\\
 = &\left(\frac{\partial}{\partial z}U(z,w) W'_{\text{out}}(z)\right.
\\&\;\;\; 
\left.+\frac{\partial}{\partial w}U(z,w) W'_{\text{in}}(z)\right)\Big|_{z=1,w=1}+1.
 \end{aligned}
 \end{equation}
Further on, differentiating  the equations \eqref{eq:Wout} and \eqref{eq:Win}, applying the definitions \eqref{eq:Ux}, \eqref{eq:Uy},  and evaluating at point $z=1$ yields the explicit expressions for $W'_{\text{in}}(1),\;W'_{\text{out}}(1),$
$$W'_{\text{out}}(1)=\frac{N_1}{A}, \;  W'_{\text{in}}(1)=\frac{N_2}{A},$$
and consequently,
\begin{equation*}
\tag{\ref{eq:derW}$'$}
\label{eq:new}
W'(1)=\frac{\mu (N_1+N_2)}{A}+1,
\end{equation*}
where 
\begin{multline}
A:=\left(2\mu\frac{\partial^2 }{\partial z \partial w  }U(z,w)-\left(\frac{\partial^2 }{\partial z \partial w  }U(z,w)\right)^2\right.
\\ 
 \left. +\left( \frac{\partial^2 }{ \partial z^2  }U(z,w) \right) \frac{\partial^2 }{\partial w^2  }U(z,w)\right)\Big|_{z=1,w=1}
\end{multline}
\begin{multline}
N_1:=\mu \frac{\partial}{\partial w}U(z,w) -\left( \frac{\partial^2}{ \partial w\partial z}U(z,w)\right)  \frac{\partial}{\partial w}U(z,w)
\\ 
+\left(\frac{\partial^2}{\partial w^2}U(z,w) \right) 
\frac{\partial}{\partial z }U(z,w)\Big|_{z,w=1}
\end{multline}
\begin{multline}N_2:=\mu \frac{\partial}{\partial z}U(z,w)-\left(\frac{\partial^2}{ \partial w \partial z}U(z,w) \right) \frac{\partial}{\partial z}U(z,w) 
\\ 
+\left( \frac{\partial^2}{\partial z^2}U(z,w)\right)  \frac{\partial}{\partial w }U(z,w)\Big|_{z,w=1}
\end{multline}
Now, by looking at the structure of \eqref{eq:new}, we see that this expression diverges when $A \to 0.$ This is the point that marks the phase transition as it implies a singularity in the average component size.
Definitions of the moments \eqref{eq:MoM} allow us to rewrite $A$ in a shorter form: 
 the directed random graph contains the giant \emph{weak} component iff 
\begin{equation}
\label{eq:criterion}
A=2\mu \mu_{11}- \mu\mu_{02}   - \mu \mu_{20} + \mu_{02}   \mu_{20} - \mu_{11}^2>0.
\end{equation}
It is interesting to compare this result to similar findings for other types of giant components, studied elsewhere, see Table~\ref{tab:criteria}. 
For instance, in undirected random graphs there is only one notion for connected component, and the \emph{giant component}\cite{molloy1995} exists iff
\begin{equation}
\label{eq:criterionMolloy}
\mu_{2}-2\mu_1 >0,
\end{equation}
where $\mu_1,\mu_2$ are the first two moments of the degree distribution. If the degree distribution is simply a translated discrete delta function, $d(l)=\delta(l-k)$, thus having $\mu_1=k,\;\mu_2=k^2,$ the criterion \eqref{eq:criterionMolloy} degenerates to $k\geq3.$  From the perspective of percolation theory this means that a regular Bethe lattice\cite{fisher1961} admits unbounded clusters of infinite size only if the corresponding coordinate number, $\sigma=k-1\geq2$. Less trivial degree distributions generalize this expression to irregular Bethe lattices. In order to derive the critical probability, one needs a dynamic process that assigns a specific degree distribution to a measure of progress $c$ (or in the case of percolation on Bethe lattices, probability $p$). Such process will be discussed in Section III.

In directed graphs there are three types of connected components: in-components, out-components, and weak components. 
The giant \emph{in-component}\cite{newman2001} exists iff, 
\begin{equation}
\label{eq:criterionNewman}
\mu_{11}-\mu >0.
\end{equation}
This inequality is stronger then the criterion \eqref{eq:criterion}, meaning that the existence of the giant in-component is also sufficient for giant weak component to exist.   Furthermore, the criterion for existence of giant \emph{out-component} is identical to \eqref{eq:criterionNewman}.

For a given directed degree distribution $u(n,k)$, we may associate a one-dimensional degree distribution by disregarding the direction of edges, 
$d(l)=\sum_{n+k=l} u(n,k)$. In this case one may apply \eqref{eq:criterionMolloy} to find out if the giant component exists in the induced undirected graph.\cite{dorogovtsev2001} When expressed in terms of moments of the bivariate distribution $u(n,k)$ this criterion shapes as
$2 \mu_{11} + \mu_{02} + \mu_{20} - 4 \mu > 0.$
The  criterion, however, should not be interpreted as the existence criterion for the giant weak component as it refers to a different topology.
\begin{table*}
\begin{tabular}{ccc}
 Type  & Criteria & Reference\\
\hline
\parbox[l]{5cm}{\begin{flushleft}  undirected graphs, percolation on Bethe lattices: \emph{giant component}\\ \end{flushleft}}& $\mu_2 - 2 \mu_1>0$ &
\parbox[l]{5cm}{ Molloy \& Reed\cite{molloy1995}
\\
Fisher and Essam\cite{fisher1961}}
\\
\parbox[l]{5cm}{ \begin{flushleft}  directed graphs: \\ \emph{giant in-component}, \\ \emph{giant out-component } \end{flushleft} }     & $\mu_{11}-\mu>0$ & \parbox[l]{4cm}{Newman, Strogatz,\\ Watts \cite{newman2001} }\\
\parbox[l]{5cm}{\begin{flushleft}  directed graph:\\
\emph{giant weak component }\end{flushleft} }  &\parbox[l]{4cm}{ \begin{align*}2\mu \mu_{11}- \mu\mu_{02}   - \mu \mu_{20} &\\+ \mu_{02}   \mu_{20} - \mu_{11}^2&>0\end{align*}} & This work \\
\end{tabular}
\caption{Existence criteria  for various types of giant components in directed and undirected graphs as a function of degree-distribution moments.}
\label{tab:criteria}
\end{table*}

\section{Evolving directed graphs with arbitrary bounded degrees}
\begin{figure}[htbp]
\begin{center}
\includegraphics[width=0.3\textwidth]{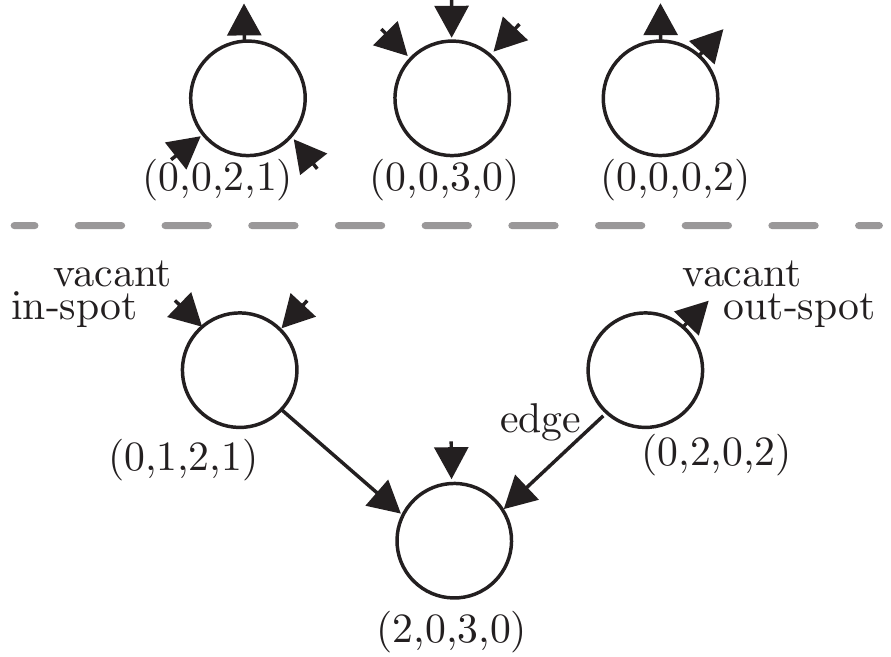}
\caption{Each vertex of the evolving directed graph is viewed as carrying two types of edges (in/out) and two types of vacant spots (also in/out).
Initial number of vacant spots sets up a limit on the maximum number of incident edges for each vertex.
The random process converts a pair of an in-spot and an out-spot belonging to different vertices into a directed edge.}
\label{fig:spots_and_edges}
\end{center}
\end{figure}

In this section we construct a time-continuous random process for evolution of the directed random graph. A specific feature of this process is that the in-/out-degree of each vertex is bounded according to a priori specified distribution. The state of each vertex is described by vector $(n,k,n_{\text{max}},k_{\text{max}})$, where $n$ counts in-edges, $k$ -- out-edges, and $n_{\text{max}},\, k_{\text{max}}$ are bounds on the maximum numbers for edges of each type. The bounds are not the same for different vertices, and initially, when no in-/out-edges are present, the whole system is characterized only by distribution of bounds $P(n_{\text{max}},k_{\text{max}}): \mathbb{N}^2_0\rightarrow \mathbb{R}^+$.
As the time, $t,$ progresses continuously, the vertex states are evolving according to the mechanism, 
\begin{multline}\label{eq:process}
(n_1,k_1,n_{\text{max},1},k_{\text{max},1}) + (n_2,k_2,n_{\text{max},2},k_{\text{max},2}) 
\rightarrow\\ (n_1+1,k_1,n_{\text{max},1},k_{\text{max},1}) + (n_2,k_2+1,n_{\text{max},2},k_{\text{max},2}) 
\end{multline}
where the rate is $\tau(n_{\text{max},1}-n_1) (k_{\text{max},2}-k_2).$ Here $\tau$ is a rate constant that with no loss of generality may be considered to be unity.
The difference between the edge bound and the actual number of edges,  $n_{\text{max},1}-n_1$ (or $k_{\text{max},2}-k_2$) refers to the finite capacity of a vertex to receive a new edge. In this process, not every pair of vertices have an equal probability to become connected but the vertices that have greater capacity are preferred. As an alternative notation for a vertex state, one may thus speak of vertices with $n-n_{\text{max}}$ vacant spots for in-edges, and $k_{\text{max}}-k$ vacant spots for out-edges, as illustrated in Fig.~\ref{fig:spots_and_edges}.

We will now study how the evolution process \eqref{eq:process} affects the degree distribution.
\subsection{Distribution for degrees and degree bounds}
As a matter of convention, the degree distribution for a directed graph is a two-dimensional function counting in-degree $n$ and  out-degree $k$. In the current case, we employ a degree distribution with two extra dimensions to account for bounds on in-degree $(n_{\text{max}})$  and out-degree $(k_{\text{max}})$.
At any point in time $t \geq 0,$ probability measure $u(n,k,n_{\text{max}},k_{\text{max}},t) \in \mathbb{N}_0^2 \times \mathbb{N}^2\supset \Omega   \rightarrow \mathbb{R}^+$ denotes the probability of randomly selecting a vertex with in-degree $n,$ out-degree $k$, and in/out-degree bounds $n_{\text{max}},k_{\text{max}}>0.$ 
Thus the values of state vector $(n,k,n_{\text{max}},k_{\text{max}})\in \Omega,$  always satisfy $0 \leq n \leq n_{\text{max}},$ and $0 \leq k\leq k_{\text{max}}.$
As the process evolves, degrees $n,k$ of each node depart from 0 and increase, while the degrees bounds $n_{\text{max}},k_{\text{max}}$ remain fixed. For this reason moments of the probability measure $u$ over $n_{\text{max}},k_{\text{max}}$ are time-independent,
\begin{equation}
\label{eq:mu0011}
\nu_{ij} = \sum_\Omega n_{\text{max}}^iy_{\text{max}}^j \,u(n,k,n_{\text{max}},k_{\text{max}},t),
\end{equation}
while the moments over $n,k$ are functions of time,
\begin{equation}
\label{eq:mu1100}
\mu_{ij}(t)=\sum_\Omega n^iy^j \,u(n,k,n_{\text{max}},k_{\text{max}},t),\\
\end{equation}
The first moments have a clear interpretation: $\nu_{10}$,\,$\nu_{01}$  denote total numbers of vacant spots for in- and out-edges present initially, $\mu_{10}$,\,$\mu_{01}$ denotes total number of in-/out-edges present in the graph at the current instance of time. 
The total probability is conserved and equal to the total probability of the initial distribution,
 \begin{multline}
\label{eq:0moment}
\sum_\Omega u(n,k,n_{\text{max}},k_{\text{max}},t) =\\ \sum_{n_{\text{max}},k_{\text{max}}} P(n_{\text{max}},k_{\text{max}})=1,\; t \geq 0.
\end{multline}
In the directed graph, the total numbers of in-edges and out-edges coincide, which in terms of our notation casts out as an additional constrain on the degree distribution, 
$$\mu(t)=\mu_{10}(t)=\mu_{01}(t),\; t>0.$$
Note, a similar equality for the degree bounds, generally speaking, does not hold ($\nu_{10} \neq \nu_{01}$) since one is free to choose initial conditions arbitrary. Yet, in the partial case when $\nu_{10} = \nu_{01},$ the initial distribution $P(n,k)$ also defines a valid graph topology that is an irregular Bethe lattice. In this case, the process is equivalent to a mean-field percolation process on this lattice. In percolation theory, As a matter of convention, instead of time as a measure of the progress one employs the edge occupancy probability, $p=\frac{\mu(t)}{\nu_{10}(t)}=\frac{\mu(t)}{\nu_{01}(t)}.$ 

The total number of vacant in-spots in the whole system is given by $\nu_{10} -\mu(t),$ and $\nu_{01} -\mu(t)$ refers to out-spots.  When a vertex $(n,k,n_{\text{max}},k_{\text{max}})$ receives a new out-edge, the choice is made between $n_{\text{max}}-n$ vacant out-spots on the vertex and 
 $\nu_{01} -\mu(t)$ vacant in-spots in the whole system, thus the rate $(n_{\text{max}}-n)(\nu_{01} -\mu(t))$.
 Similar considerations are made for the placement of an in-edge, and the dynamics for the degree distribution shapes as the following master equation, 
\begin{equation}
\label{eq:balance}
\begin{aligned}
\frac{\partial}{\partial t} u(n,k,t)  =&
  (n_{\text{max}} - n + 1) (\nu_{01} - \mu(t)) u(n-1,k,t) \\
& +(k_{\text{max}} - k + 1) (\nu_{10} - \mu(t)) u(n,k-1,t)\\
& -\Big(\,(n_{\text{max}} - n) (\nu_{01} - \mu(t)) \\
 &\;\;\;\,\;\;\;\;+(k_{\text{max}} - k) (\nu_{10} - \mu(t) )\,\Big)u(n,k,t).
\end{aligned}
\end{equation}
$$u(n,k,n_{\text{max}},k_{\text{max}},0)  =\delta(n)\delta(k) P(n_{\text{max}},k_{\text{max}})$$
Here, discrete delta functions in the initial conditions, $\delta(n)\delta(k),$ refer to the fact that vertices have no edges initially. 
Further in the text, where it leads to no confuse, we will drop $n_{\text{max}},k_{\text{max}}$ dimensions referring to $u(n,m,t)$ for the sake of a shorter notation.
As a routine to solve \eqref{eq:balance} we take the following steps: first, the differential-difference equation \eqref{eq:balance}  will be transformed to a non-linear PDE by the generating function transform; 
then, we derive and solve an ODE for $\mu(t)$ that also eliminates the nonlinearity;
finally, the linear PDE is solved and the solution is transformed back to the domain of discrete functions.
Below, these steps are elucidated in more details.

We act on left and right hand sides of the balance equation \eqref{eq:balance} with the bivariate generation function transform \eqref{eq:transform} in dimensions $n,k$.
Thus in the generating-function domain \eqref{eq:balance} becomes,
\begin{equation}\label{eq:non-linearPDE}
\begin{aligned}
\frac{\partial }{\partial  t}&U(z,w,t) = \\
&
 \Big( n_{\text{max}}( z-1) \nu_{01}+ k_{\text{max}}(w-1) \nu_{10} \\& \;\;+(n_{\text{max}}+k_{\text{max}}-  k_{\text{max}} z- k_{\text{max}} w) \mu(t)  \Big) U(z,w,t)\\
&-z ( z-1 ) (\nu_{01}-\mu(t)) \frac{\partial}{\partial z}U(z,w,t) \\
&-(w-1) w (\nu_{10}-\mu(t)) \frac{\partial}{\partial w}U(z,w,t).
\end{aligned}
\end{equation}
\begin{multline*}
U(n,k,n_{\text{max}},k_{\text{max}},t)|_{t=0}  = P(n_{\text{max}},k_{\text{max}}),\; 
\\ (n,k,n_{\text{max}},k_{\text{max}})\in\Omega.
\end{multline*}
This PDE is not linear as the unknown function, $U(z,w,t),$ is also used in the definition of $\mu(t)$. We can bring \eqref{eq:non-linearPDE} to a simpler form by first resolving an expression for $\mu(t)$.
Transform \eqref{eq:transform} maps weighted distributions to partial derivatives of the corresponding generating functions (\textit{e.g.} $n u(n,k,t) \rightarrow z \frac{\partial}{\partial z} U(z,w,t)$), and  the sum over the whole domain to the value of the generation function  at point 1, (\textit{e.g.} $\sum_{n,k} u(n,k,t)\rightarrow U(z,w,t)|_{z=1}).$
Thus applying operator $z\frac{\partial}{\partial z}\cdot |_{z=1}$ to the both sides of \eqref{eq:non-linearPDE} yields an ODE for the first moment:
\begin{equation}\label{eq:ode1000}
\left\{
\begin{aligned}
\mu'(t)=&\big(\nu_{01} - \mu(t)\big) \big(\nu_{10} -  \mu(t)\big),\\
\mu(0)=&0.
\end{aligned}
\right.
\end{equation}
The the solution of the differential equation \eqref{eq:ode1000} reads,
\begin{equation}\label{eq:1000}
\mu(t)=\nu_{01}-\frac{\nu_{01} (\nu_{01}-\nu_{10})}{\nu_{01}-\nu_{10}\, \text{e}^{t (\nu_{10}-\nu_{01})}}.
\end{equation}
Having an expression for $\mu(t)$ that contains only $t$ and constants allows us to separate the variables in \eqref{eq:non-linearPDE}. Assuming $$U(z,w,n_{\text{max}},k_{\text{max}},t)=f(z,n_{\text{max}},k_{\text{max}},t)g(w,n_{\text{max}},k_{\text{max}},t)$$ and dropping $n_{\text{max}},k_{\text{max}}$ dimensions in the shorthand notation, we obtain
\begin{equation}
\label{eq:fg}
\begin{aligned}
&\frac{A_1(z) f(z,t)+B_1(z)\frac{\partial}{\partial z} f(z,t)-\frac{\partial}{\partial t}f(z,t)}{f(z,t)}=\alpha(t)
\\
&\frac{A_2(w) g(w,t)+B_2(w) \frac{\partial}{\partial w}g(w,t)-\frac{\partial}{\partial t}g(w,t)}{g(w,t)}=-\alpha(t)
\end{aligned}
\end{equation}
where $\alpha(t)$ does not depend on $z,w,$ and the rest of the coefficients are as follows,
$$
\begin{aligned}
A_1(z)&=\big(\nu_{01} -  \mu(t)\big)n_{\text{max}} (z-1), \\
A_2(w)&= \big(\nu_{01} -  \mu(t)\big)k_{\text{max}} (w-1), \\
B_1(z)&= \big(\nu_{01} -  \mu(t)\big)z(z-1), \\
B_2(w)&= \big(\nu_{01} -  \mu(t)\big) w(w-1). \\
\end{aligned}
$$
Additionally, a solution of PDE \eqref{eq:fg} have to satisfy the total probability conservation \eqref{eq:0moment}. To ensure that, we apply operators $\sum_{n_{\text{max}},k_{\text{max}}}\cdot |_{z=1}$, $\sum_{n_{\text{max}},k_{\text{max}}}\cdot |_{w=1}$ to both parts of \eqref{eq:fg}:
$$
\begin{aligned}
-\frac{\partial}{\partial t}\sum_{n_{\text{max}},k_{\text{max}}}f(1,t)=\alpha(t) \sum_{n_{\text{max}},k_{\text{max}}}f(1,t).\\
\frac{\partial}{\partial t}\sum_{n_{\text{max}},k_{\text{max}}}g(1,t)=\alpha(t)\sum_{n_{\text{max}},k_{\text{max}}}g(1,t).
\end{aligned}
 $$
 From here, it becomes obvious that the only $\alpha(t)$ that admits the total probability conservation  for non-zero $f(z,t),g(w,t)$ is $\alpha(t)\equiv0.$
Therefore, the PDEs introduced in \eqref{eq:fg} simplify to
\begin{equation}\left\{
\begin{aligned}
\frac{\partial}{\partial t}f(z,t)=&A_1(z) f(z,t)+B_1(z)\frac{\partial}{\partial z} f(z,t),\\
f(z,t)|_{t=0} = &P_1(n_{\text{max}});
\end{aligned}\right.
\end{equation}
\begin{equation}
\left\{
\begin{aligned}
\frac{\partial}{\partial t}g(w,t)=&A_2(w) g(w,t)+B_2(w) \frac{\partial}{\partial w}g(w,t),\\
g(w,t)|_{t=0} =& P_2(k_{\text{max}}).
\end{aligned}\right.
\end{equation}
and lead to the following solutions,
\begin{equation}
\label{eq:Usol}
\begin{aligned}
f(z,t)=&\left(1+\frac{(z-1)}{\nu_{10} } \mu(t)\right)^{n_{\text{max}}}P_1(n_{\text{max}}),\\
g(w,t)=&\left(1+\frac{(w-1) }{\nu_{01}}\mu(t)\right)^{k_{\text{max}}}P_2(k_{\text{max}}),
\end{aligned}
\end{equation}
where $P(n_{\text{max}},k_{\text{max}})=P_1(n_{\text{max}})P_2(k_{\text{max}}).$
Having expressions for $f(z,t),g(w,t)$ permits a straightforward asymptotical analysis for $u(n,k,t)$ when $t$ approaches infinity. Depending on the relation between parameters $\nu_{01},\nu_{10}$ three modes emerges.\\
1)Equal maximum numbers of in- and out-edges, $\nu_{01}=\nu_{10}.$ In this case, 
$$\lim_{\nu_{01}\to\nu_{10}}\mu(t) = \frac{\nu_{10}^2 t}{(1 + \nu_{10} t)^2},$$
and $$\lim_{t\to \infty}U(z,w,t) = z^{n_{\text{max}}}w^{k_{\text{max}}}P(n_{\text{max}},k_{\text{max}}).$$
This expression generates a discrete delta function translated to position $(n_{\text{max}},k_{\text{max}})$,
\begin{multline}
\label{eq:delta}
u(n,k,n_{\text{max}},k_{\text{max}},t)=\\ \delta(n-n_{\text{max}})\delta(k-k_{\text{max}})P(n_{\text{max}},k_{\text{max}}).
\end{multline}
2)The maximum number of in-edges exceeds the number of out-edges, $\nu_{10}>\nu_{01}.$ In this case, the distribution cannot evolve into a single delta function at infinite time, since there will always be vertices with unused spots for an in-edge. In fact, the evolution of the system will stop when, $\nu_{01}=\nu_{10}.$
As $t\to \infty$ the generation function approaches to
$$\lim_{t\to \infty}U(z,w,t) = \left( 1 +  \frac{\nu_{01}}{\nu_{10}} (z-1)\right)^{n_{\text{max}}}w^{k_{\text{max}}}P(n_{\text{max}},k_{\text{max}}),$$
which generates
\begin{multline}
\label{eq:assymptotical_x}
\lim_{t\to\infty}u(n,k,n_{\text{max}},k_{\text{max}},t)=\\ \binom{n_{\text{max}}}{n}\left(\frac{\mu(t)}{\nu_{10}}\right)^n \left(1 - \frac{\mu(t)}{\nu_{10}}\right)^{n_{\text{max}}-n} \!\!\!\!\!\!\!\delta(k-k_{\text{max}})P(n_{\text{max}},k_{\text{max}}).\end{multline}

3)The maximum number of out-edges exceeds the number of in-edges, $\nu_{01}>\nu_{10}.$ Analogously to the previous case, one obtains the limiting value for the degree distribution,
\begin{multline}
\lim_{t\to\infty}u(n,k,n_{\text{max}},k_{\text{max}},t)=\\
\binom{k_{\text{max}}}{k}\left(\frac{\mu(t)}{\nu_{01}}\right)^k \left(1 - \frac{\mu(t)}{\nu_{01}}\right)^{k_{\text{max}}-k} \!\!\!\!\!\!\!\delta(n-n_{\text{max}})P(n_{\text{max}},k_{\text{max}}).
\end{multline}
Finally, in the general case of finite time, the expression for the degree distribution is generated by \eqref{eq:Usol},
\begin{equation}
\label{eq:u}
\begin{aligned}
u(n,&k,n_{\text{max}},k_{\text{max}},t)=\\
&\binom{n_{\text{max}}}{n}\binom{k_{\text{max}}}{k}\left(\frac{\mu(t)}{\nu_{10}}\right)^n \left(1 - \frac{\mu(t)}{\nu_{10}}\right)^{n_{\text{max}}-n}\\
&\times\left(\frac{\mu(t)}{\nu_{01}}\right)^k \left(1 - \frac{\mu(t)}{\nu_{01}}\right)^{k_{\text{max}}-k}\!\!\!P(n_{\text{max}},k_{\text{max}}).
\end{aligned}
\end{equation}
\begin{figure*}[htbp]
\begin{center}
\includegraphics[width=0.95\textwidth]{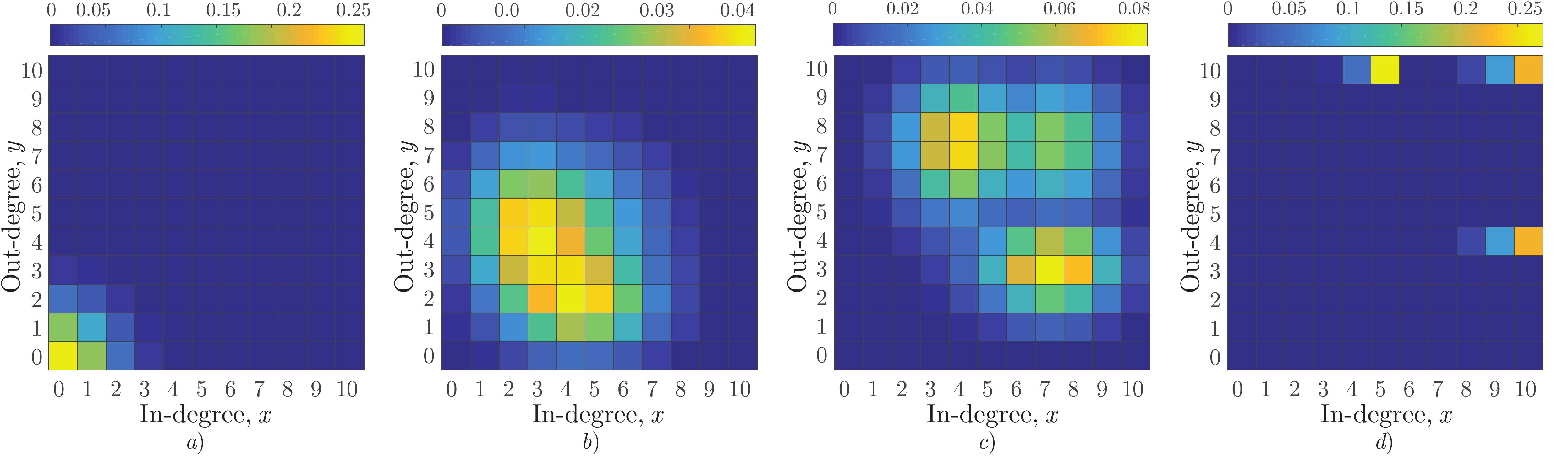}
\caption{
Evolution of the bivariate degree distribution $u(n,k,t)$ for the random graph model with bounded degrees.
In this example, the initial distribution of bounds $P(n_{\text{max}},k_{\text{max}})$ vanishes everywhere except for the points: $P(10,10)=\frac{1}{3},\;P(5,10)=\frac{1}{3},\;P(10,4)=\frac{1}{3}.$ Since the total number of in-spots exceeds the total number of out-spots,  only $k-$marginal of the distribution approaches a linear combination of delta functions in the time limit, $t\to \infty.$ 
The time snapshots are obtained for the following values of time: a) $t=0.01$, b) $t=0.1$, c) $t=1$, d) $t\to\infty.$
}
\label{fig:cubes}
\end{center}
\end{figure*}
Although four-dimensional  distribution $u(n,k,n_{\text{max}},k_{\text{max}})$ has a relatively simple expression when viewed for a specific class of vertices $n_{\text{max}},k_{\text{max}},$ the most useful output of this model is the two-variate degree distribution $d(n,k,t)=\sum_{n_{\text{max}},k_{\text{max}}}u(n,k,n_{\text{max}},k_{\text{max}},t).$ Distribution $d(n,k,t)$ may exhibit a 'non-trivial' interplay of peaks for certain initial distributions. Fig.~\ref{fig:cubes} illustrates evolution of the degree distribution for a sample system.  The initial distribution 
of spots, $P(n_{\text{max}},k_{\text{max}}),$ is non-zero only in three points $P(10,10), P(5,10), P(10,4).$  Naturally, evolution of $d(n,k)$ starts with all probability density located at point $(0,0)$ at $t=0$. In the intermediate time stages the distribution becomes broad, so that $P(n,k)>0,\;n,k\leq10$. Asymptotically, $d(n,k,t)$ converges to a steady state at $t\to \infty$. Since there are more possibilities for in-edges than out-edges, $\nu_{10}>\nu_{01}$, the steady-state degree distribution is of the type \eqref{eq:assymptotical_x}.

\subsection{Phase transition}
Now, when we have an explicit expression for the degree distribution, it is possible to apply the existence criterion  \eqref{eq:criterion} to this expression and in this way find the critical parameters.
Let $c_n(t)$ and $c_k(t)$ denote the fraction of in-spots and out-spots that were converted into in-edges,
 \begin{equation}
\begin{aligned}
 c_n(t)= \frac{\mu(t)}{\nu_{10}},\\
 c_k(t) = \frac{\mu(t)}{\nu_{01}}.
 \end{aligned}
 \end{equation}
Both, $c_n(t)$ and $c_k(t)$ are non-negative; their upper bounds, however, depend on the initial conditions. 
Namely, 
 \begin{equation}
 \label{eq:bounds}
 \begin{aligned}
 \sup_{t>0}  c_n(t)= 
 \begin{cases}
 1,                                               &\nu_{01} \geq \nu_{10},\\
 \frac{\nu_{01}}{\nu_{10}}, & \nu_{01}<\nu_{10};
 \end{cases}\\
  \sup_{t>0} c_k(t)= 
 \begin{cases}
 \frac{\nu_{10}}{\nu_{01}}, &\nu_{01} \geq \nu_{10},\\
 1, & \nu_{01}<\nu_{10}.
 \end{cases}
 \end{aligned}
 \end{equation}
In order to retrieve explicit expressions for moments that appear in \eqref{eq:criterion}, we act on the solution \eqref{eq:u} with differential operators 
$\left( \frac{\partial^2 }{\partial z^2} - \frac{\partial }{\partial z}\right), $
$ \left( \frac{\partial^2 }{\partial w^2} - \frac{\partial }{\partial w}\right),$ and
$ \left( \frac{\partial^2 }{\partial z\partial w}\right) $
to  correspondingly obtain:
 \begin{equation}
 \label{eq:mom2}
 \begin{aligned}
\mu_{02}(t) =& c_k(t) \nu_{01} - c_k(t)^2 \nu_{01} + c_k(t)^2 \nu_{02},\\
\mu_{20}(t) = &c_n(t) \nu_{10} - c_n(t)^2 \nu_{10} + c_n(t)^2 \nu_{20},\\
\mu_{11}(t) = &c_n(t) c_k(t) \nu_{11}.\\
 \end{aligned}
 \end{equation}
Plugging  \eqref{eq:mom2} into \eqref{eq:criterion} and realizing that $c_k(t)=\frac{\nu_{10}}{\nu_{01}}c_n(t)$ yields a criterion for existence of the giant weak component as a quadratic function of $c_n(t)$ with the coefficients involving exclusively initial moments, 
 \begin{equation}
 \label{eq:quad}
 a \,c_n(t)^2 +b \,c_n(t) + c>0,
  \end{equation}
where
$$
 \begin{aligned}
a=& \nu_{01} \nu_{10} - \nu_{02} \nu_{10} - \nu_{11}^2 -  \nu_{01} \nu_{20} + \nu_{02} \nu_{20},\\
b= &2  \nu_{01} \nu_{11},\; c=-\nu_{01}^2.\\
 \end{aligned}
$$
The both roots of \eqref{eq:quad} are real, and the smallest root is always negative. Therefore, the inequality is satisfied when,
  \begin{equation}
  \label{eq:weakgiantx}
 c_n(t) >c_{n,\text{critical}}= \frac{\nu_{01}}{ \nu_{11} +\sqrt{(\nu_{02}-\nu_{01}) (\nu_{20}-\nu_{10})}}
   \end{equation}
or   alternatively,
  \begin{equation}
  \label{eq:weakgianty}
 c_k(t) >c_{k,\text{critical}}= \frac{\nu_{10}}{ \nu_{11} +\sqrt{(\nu_{02}-\nu_{01}) (\nu_{20}-\nu_{10})}}
   \end{equation}
Equations \eqref{eq:weakgiantx} and \eqref{eq:weakgianty} express the main result of this section: the phase-transition point in terms of a monotone function of time, $c_n(t).$ One may easily transit to the actual time, $t,$ by evaluating 
$$t = \frac{\log \left( \frac{(1 - c_n) \nu_{01}}{\nu_{01} - c_n \nu_{10}}\right)}{\nu_{10} - \nu_{01}},\text{ or }t = \frac{\log \left( \frac{ \nu_{10} -c_k \nu_{01}   }{(1 - c_k) \nu_{10} }  \right)}{\nu_{10} - \nu_{01}}.$$
In the context of edge percolation on Bethe lattices (that this problem degenerates to when $\nu_{01}=\nu_{10}$) the both critical values \eqref{eq:weakgiantx} and \eqref{eq:weakgianty} coincide and the critical probability is given by $p_c=c_{n,\text{critical}}=c_{k,\text{critical}}.$

Finally, let us turn to the next question: does a specific initial distribution, $P(n_{\text{max}},k_{\text{max}}),$ yield the giant weak component in finite time or the phase transition never happens?
To answer this question it is enough to check if the giant component exists when $c_n$ or $c_k$ approach  their upper bounds \eqref{eq:bounds}.
Evaluating the moment expressions  \eqref{eq:mom2} at the upper bounds for $c_n,c_k,$ \eqref{eq:bounds} and substituting the expressions into the phase-transition criterion \eqref{eq:criterion} yields the desired condition: initial distribution $P(n_{\text{max}},k_{\text{max}})$ admits the phase transition in finite time iff at least one of the following conditions is true:
\begin{equation}
 \label{eq:criterion2}
 \begin{aligned}
A_1:=(\nu_{02} - \nu_{01}) (\nu_{20} - \nu_{10}) - (\nu_{11}& - \nu_{01})^2>0, \\
&\text{and }   \nu_{01} \geq \nu_{10},\\
\text{or     \hspace{4.8cm}} &\\
A_2:=(\nu_{02} - \nu_{01}) (\nu_{20} - \nu_{10}) - (\nu_{11} &- \nu_{10})^2>0,  \\
&\text{and }   \nu_{01} \leq  \nu_{10}.
\end{aligned}
 \end{equation}
 Furthermore, the asymptotic phase transition occurs at $t\to\infty$ iff the inequalities in \eqref{eq:criterion2}  are replaced by equalities (\textit{i.e.} $A_1=0,\; A_2=0$).
 
 When $\nu_{01}=\nu_{10}$ (\textit{i.e.} equal number of in-spots and out-spots are present initially) the both inequalities from \eqref{eq:criterion2}  degenerate to 
$$
 2 \nu_{10} \nu_{11} -\nu_{02} \nu_{10} + \nu_{02} \nu_{20} - \nu_{10} \nu_{20} - \nu_{11}^2  >0.
$$
One may see that this condition is identical to \eqref{eq:criterion}, this similarity is not a coincidence. When equal numbers of vacant spots for in- and out edges are used, all spots will be converted into edges at infinite time. Furthermore, the degree distribution $u(n,k,t)$ degenerates to delta function, as was given in \eqref{eq:delta}, and consequently its moments approach the moments of the distribution of degree bounds, $P(n_{\text{max}},k_{\text{max}}),$ \textit{i.e.} $\mu(t) \to \nu_{10},\;\mu_{20}(t) \to \nu_{20},\;\mu_{11}(t) \to \nu_{11},$ \textit{etc}.
\begin{figure*}[htbp]
\begin{center}
\includegraphics[width=0.9\textwidth]{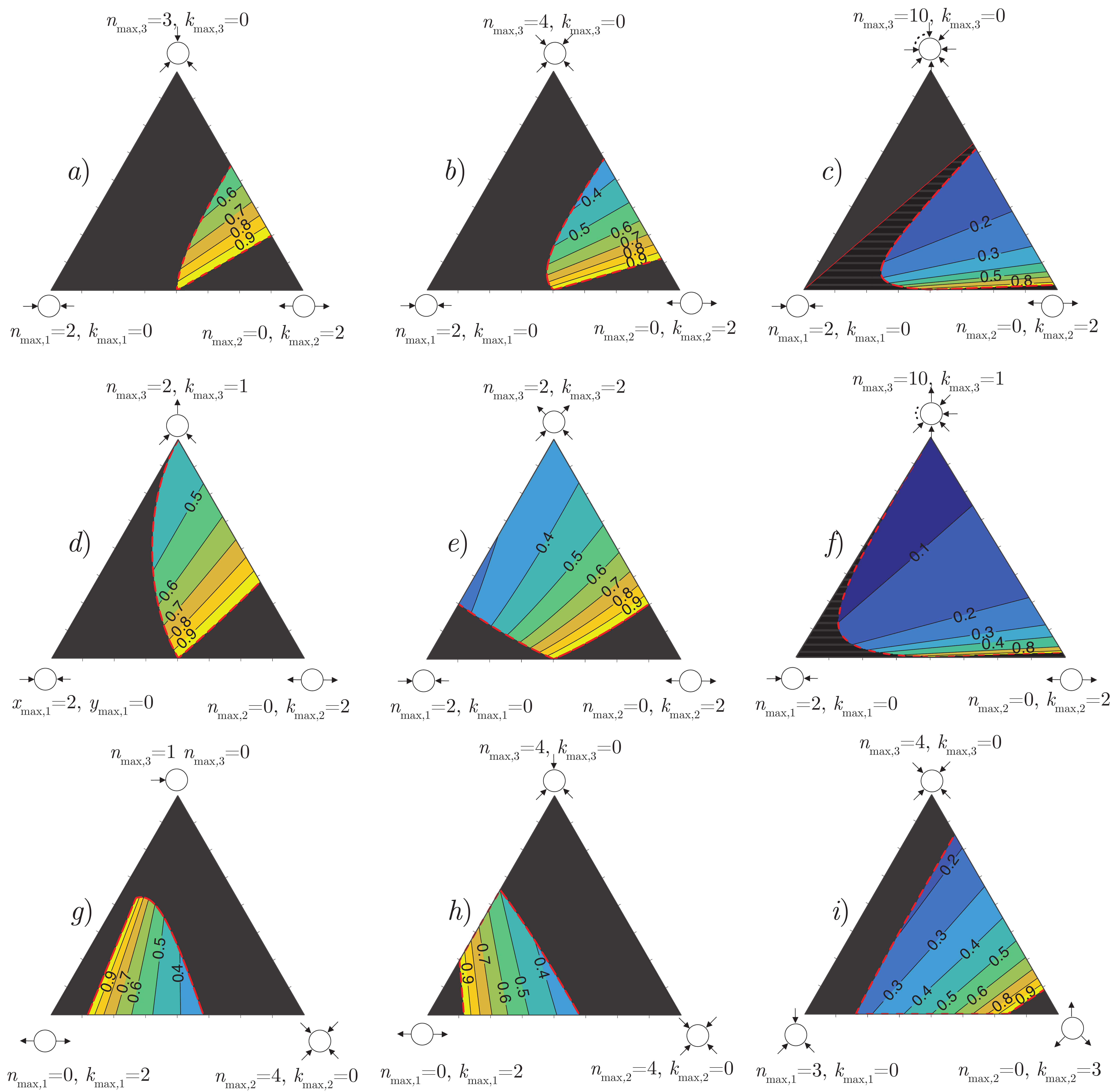}
\caption{Barycentric plots for the weak-component phase-transition time in terms of $c_n,$ as obtained for the random graph model with bounded degrees.  
 Panels ($a-i$) regard various initial distributions that are non-zero at three points $n_{\text{max},i},k_{\text{max},i},\;i=1,2,3$.  Cases $a,b$ are covered by Flory-Stockmayer theory. The black area contains configurations that do not admit phase transition, the color (shaded) area contains configurations with finite phase transition, the red dashed line contains configurations that admit the phase transition asymptotically at infinite time.}
\label{fig:bar}
\end{center}
\end{figure*}
Generally speaking, there are many possible configurations of initial distribution $P(n_{\text{max}},k_{\text{max}}).$ However, when the initial distribution is non-zero only in three points $P(n_{\text{max},i},k_{\text{max},i})=c_i,\;i=1,2,3,$ it is convenient to visualize the result of \eqref{eq:weakgiantx} for all $f_1+f_2+f_3=1$ with a plot in the barycentric coordinates (Fig.~\ref{fig:bar}). 
Since $c_i\geq0,$ the plot is contained within a triangle; the $i^{\text{th}}$ vertex of the triangle is associated with $c_i=1$ (so that the other two values are zero), and the points inside refer to all $c_i$ being non-zero.
The color (shaded) area in the panels of Fig.~\ref{fig:bar} denotes the phase-transition point in terms of $c_n$, the black area denotes configurations that lead to no phase transition and the red (dashed) line contains configurations that admit the phase transition asymptotically.

We will now focus on some qualitative properties of the cases presented in Fig.~\ref{fig:bar} for an illustrative purpose.
One may observe that when the degree bounds restrict a vertex to be a sink (only in-edges) or a source (only out-edges), the configurations that admit phase transition occupy only the area close to the centre of the triangle (cases $a-c$, $g-i$) but not close to triangle's vertices. This means that only a combination of sinks and sources yields a system with the phase transition. 
Asymptotically, when sinks and sources with maximum of two edges are combined with sinks of infinitely high degree, the phase space slits into two regions as shown in  Fig.~\ref{fig:bar}$c$: no phase transition (black area) and immediate phase transition (hatched area).
Vertices that have at least two edges of one kind and one of the other kind (\textit{i.e.} 2 out- one in-edge, or two in- one out-edge) can form the giant weak component alone. In this case, as opposed to the non-directional case, large proportion of sinks (or sources) may postpone the phase transition, see panels $d,e,f,h.$ 
Finally, if a vertex is allowed to have one edge at most, this vertex will significantly postpone emergence of the giant component or prevent it completely (compare panels $e$ and $g$).

\subsection{Relation to Flory-Stockmayer gelation theory}
The results on the phase transition, as presented in the previous section, constitute a generalization  for Flory-Stockmayer gelation theory (FSGT).
 FSGT was developed by Flory\cite{flory1941} and Stockmayer\cite{Stockmayer1944} by means of probabilistic and kinetic arguments, respectively. 
 Later, the kinetic view on the theory was advanced by  Ziff \cite{ziff1980}. 
 FSGT predicts when step-polymerization of multifunctional monomers yields an infinite structure - the \textit{gel}. One of the limitations of the theory is that only three species of monomers are present in the mixture: two species of linear and one species of branched units.
Here by taking the random graph point of view, we demonstrate how the present results generalize applicability of FSGT to an arbitrary number of monomer species  with no constrains on their functionalities.

Flory and Stockmayer considered a polymerization model where a chemical bond may appear between a pair of reactive groups of two types, $A$  and $B$. The pair consisting of one $A-$group and one $B-$group may receive a chemical bond with equal probability, but reactions between $A$  and $A$ or $B$ and $B$ are forbidden. 
The reactive groups are carried on monomers. There are three types of monomers in the system: linear monomers with two $A$ groups, linear monomers with two $B$ groups, branched monomers with $n$ groups of type $A.$  In its essence, this model is a directed random graph model of the type \eqref{eq:process}, where a vertex resembles a monomer,  an $A-$group  -- the in-spot, $B-$group -- the out-spot, and  a chemical bond resembles a directed edge, $B \to A$. 
The initial conditions are restricted to $P(n_{\text{max}},k_{\text{max}})=0$ except for $P(1,0)=f_1$ (linear unit, $A-$monomer), $P(0,1)=f_2$ (linear unit, $B-$monomer), $P(n,0)=f_3$ (branched unit); $f_1+f_2+f_3=1$. Subsequently, the expressions for the moments of  $P(n_{\text{max}},k_{\text{max}})$  are
\begin{equation}
\label{eq:momFlory}
\begin{aligned}
\nu_{10} &= 2 f_1+n f_3, &\nu_{01} &= 2 f_2,\\
\nu_{20} &= 4 f_1+f^2 f_3,&\nu_{02} &= 4 f_2,\;\nu_{11}= 0.\\
\end{aligned}
\end{equation} 
Plugging the moments  \eqref{eq:momFlory} into the criterion \eqref{eq:criterion2} immediately gives us the condition for finite-time
emergence of the giant weak component, \textit{i.e.} gel. Polymerization system contains gel if at least one of the statements is true:
\begin{equation}
\begin{aligned}
\text{a) }& f_2 (2 f_1 - 2 f_2 + (n^2-n )f_3  )> 0\\ 
&\text{and }  2 f_2  \geq 2 f_1+n f_3;\\
\text{ or } &\\
\text{b) }&2 f_2 (2 f_1 + (n^2-n ) f_3 ) -(2 f_1 + n f_3)^2 >0\\
& \text{and }  2 f_2  \leq2 f_1+n f_3.\\
 \end{aligned}
 \end{equation}
 
Alternatively, the phase-transition condition may be rewritten as a lower bound on $c_n,$
\begin{equation}
\label{eq:cxcrit}
\begin{aligned}
c_n>c_{n,critical}&=\sqrt{\frac{ f_2  }{  f_1 +\frac{ (n^2-n ) }{2}f_3 }}.
 \end{aligned}
 \end{equation}
Now, employing the original notation used by Flory \cite{flory1941},
\begin{equation}
\alpha_c=\frac{1}{n-1},\; \rho=\frac{n f_3}{(2 f_1 + n f_3)},\; r=\frac{2 f_1+n f_3}{2 f_2},\;p_A=c_n,\;p_B=c_k
\end{equation}
and realizing that $p_B = r p_A$ we, rewrite the condition \eqref{eq:cxcrit} as
\begin{equation}
\label{eq:flo}
p_A> \sqrt{\frac{ \alpha_c   } { r(\alpha_c + \rho - \alpha_c  \rho )}}\;\;\text{or}\;\;
p_B> \sqrt{\frac{ r \alpha_c   } { \alpha_c + \rho - \alpha_c  \rho }}. 
 \end{equation}
 Here, $p_A$ (or $p_B$) measures the progress of the  process, it is the probability that A-type functionality (B-type) has been converted into a chemical bond.
 From the perspective of percolation models, $p_A, p_B$ play a similar role to the site occupancy probability, $p$, as for instance in Bethe lattice percolation model, Ref.\cite{fisher1961}
 As an alternative to \eqref{eq:flo} one may also consider a single inequality for  $\alpha=p_A^2 r(\alpha_c + \rho - \alpha_c  \rho )=\frac{p_B^2}{r} (\alpha_c + \rho - \alpha_c  \rho ),$ gelation occurs if
\begin{equation}
\label{eq:last}
\alpha > \alpha_c.
\end{equation}
Inequality \eqref{eq:last} constitutes  the central result of Flory-Stockmayer theory.

\section{Conclusions}
The fundamental assumption we rely upon in this paper is that in all respects other than their degree distribution, the graphs are treated as entirely random.
This assumption allows us to construct a powerful toolbox that connects a local property of the random graph, namely the degree distribution, to the global properties. In this respect, the inequality \eqref{eq:criterion} crystallizes as the most generic result: it allows one to verify existence of the giant weak component by knowing only moments of the degree distribution. 
No limitations are posed on the degree distribution itself. That means that one has multiple options when  applying the criterion to a particular problem.
 One option is to find the distribution empirically from measured data, which is a much easier task than measuring the weak-component size directly, for instance, in the case of social networks or the World Wide Web structure. 
Another option is to predict the degree distribution by a computer simulation, which is the method of choice in statistical mechanics among other fields.
Finally, one may apply the criterion \eqref{eq:criterion} to a theoretical model that yields an analytic expression for the degree distribution or its moments. In the latter case, the phase-transition criterion may be reformulated in terms of the model parameters.
As an example of this path, we referred to the random graph process with bounded degrees in the second part of the paper. This model plays an important role in the soft matter physics where it is used as a prototype for step-polymerization and gel formation. Instead of a computer simulation, the expression for the degree distribution is obtained analytically.
The analytic expression is then used to find the phase-transition point for the weak component in terms of the only model parameter --  the distribution of degree bounds. 
In this way it is possible to avoid resolving the whole component size distribution when focusing on the phase transition alone.

In the context of step-polymerization, the emergence of the giant weak component signifies gel formation.  In polymer synthesis, the identification of the gel point is usually associated with Flory-Stockmayer gelation theory. We showed that Flory-Stockmayer theory can be viewed as a special case of the random graph model with bounded degrees. Furthermore, as being more general, the analytic results on the random graph model with bounded degrees naturally extend Flory-Stockmayer theory to a broader scope of cases.

\begin{acknowledgments}
This work is part of the research program Veni with project number 639.071.511, which is financed by the Netherlands Organisation for Scientific Research (NWO).
\end{acknowledgments}

\bibliographystyle{apsrev4-1}
\bibliography{literature}

\end{document}